\documentclass{amsart}
\usepackage{graphicx,amssymb,amsmath,epsfig}
\usepackage{amsmath, amscd, mathrsfs,amsfonts}
\usepackage{rotating}
\usepackage{enumerate}
\newtheorem{theorem}{Theorem}[section]
\newtheorem{proposition}[theorem]{Proposition}
\newtheorem{lemma}[theorem]{Lemma}
\newtheorem{corollary}[theorem]{Corollary}
\newtheorem{remark}[theorem]{Remark}
\newtheorem{example}[theorem]{Example}

\newtheorem{definition}[theorem]{Definition}
\newtheorem{notation}[theorem]{Notation}

\usepackage{amssymb,amsmath}
\usepackage[all]{xy}

\newcommand{\bth}{\begin{theorem}}
\newcommand{\bpr}{\begin{proposition}}
\newcommand{\epr}{\end{proposition}}
\newcommand{\bco}{\begin{corollary}}
\newcommand{\eco}{\end{corollary}}
\newcommand{\ble}{\begin{lemma}}
\newcommand{\ele}{\end{lemma}}
\newcommand{\bre}{\begin{remark}\rm}
\newcommand{\ere}{\end{remark}}
\newcommand{\bex}{\begin{example}\rm}
\newcommand{\eex}{\end{example}}
\newcommand{\bde}{\begin{definition}\rm}
\newcommand{\ede}{\end{definition}}

\def\bbr{{\mathbb R}}

\def\bbr{{\mathbb R}}

\def\brho{{\mathcal B}^{Q_r}_\rho(X)}

\begin{document}

\title{Formal Barycenter Spaces with Weights: The
Euler Characteristic}

\author{Sadok Kallel}
\address{American University of Sharjah (UAE) and Laboratoire Painlev\'e, Lille 1 (France)}
\email{skallel@aus.edu}

\date{December 2018}
\maketitle

\begin{abstract} We compute the Euler characteristic with compact supports $\chi_c$ of the formal barycenter spaces with weights of some locally compact spaces, connected or not. This reduces to the topological Euler characteristic $\chi$ when the weights of the singular points are less than one. As foresighted by Andrea Malchiodi, our formula is related to the Leray-Schauder degree for mean field equations on a compact Riemann surface obtained by C.C. Chen and C.S. Lin.
\end{abstract}

\section{Statement of the Main Result}

Given a space $X$, we will write ${\mathcal B}_k(X)$ for the space of formal barycenters of $k$ points in $X$ \cite{kk}. By construction there are inclusions ${\mathcal B}_k(X)\hookrightarrow {\mathcal B}_{k+1}(X)$ for all $k$ and we will write ${\mathcal B}(X)$ the direct limit. This is known to be a contractible space if $X$ is of the homotopy type of a CW.

Let $Q_r:=\{y_1,\ldots, y_r\}\subset X$ be a fixed finite set of ``singular points" in $X$.
We assign to every $x\in X$ a \textit{weight}
$$w(x) =\begin{cases} 1,& x\not\in Q_r\\ w_i,&x=y_i\end{cases}$$
where $w_i>0$. Let $\rho $ be any positive number and define the set
\begin{equation}\label{maindef}
\brho  = \left\{ \sum t_ix_i  \in {\mathcal B}(X)\ |  \ \sum_i w(x_i)\leq \rho\right\}
\end{equation}
This is topologized as a subspace of ${\mathcal B}(X)$.
If $w_i=1$ for all $i$, the singular points are ``invisible'' (i.e. they cease to be singular) and
${\mathcal B}_{\rho}^\emptyset(X) = {\mathcal B}_{\lfloor \rho\rfloor}(X)$, where $\lfloor \rho\rfloor$ is the greatest integer less or equal to $\rho$ (i.e. $\lfloor -\rfloor$ is the floor function).

The formal barycenter spaces with weights play nowadays a significant role in geometric analysis.
They were introduced in \cite{cm} in order to study singular Liouville equations arising in the problem of prescribing the Gaussian curvature and the appearence of conical singularities on compact Riemann surfaces under a conformal change of the metric.
The weighted barycenter spaces come with filtration terms that relate to low sublevels of a $C^1$-functional whose Euler-Lagrange equation is the Liouville type equation. The non-contractibility of the weighted barycenter spaces implies a change of topology in the sublevels from which the existence of solutions is deduced. A conjecture about the contractibility of $\brho$ is stated in the case $X=\Sigma$ is a closed Riemann surface, and this conjecture is addressed in \cite{carlotto}. The computation of the Euler characteristic that we provide in this note gives precise, albeit weaker conditions on the contractibility of $\brho$ for general $X$, connected or not, compact or not. We expect that this result, in the case when $X$ is a proper smooth (not necessarily connected) subset of a compact Riemann surface $\Sigma$, can enable one to determine the Leray-Schauder degree formula for the singular Liouville equation appearing when the prescribed curvature is sign-changing, a problem recently addressed in \cite{dlr}, extending the computation done in \cite{cl} for positive curvatures.

Throughout the paper, $\chi$ will denote the topological Euler-Poincar\'e characteristic, and
$\chi_c$ the Euler characteristic with compact support (see \S\ref{compact}).
When $X$ is connected, and there are no singular weights so that $Q_r=\emptyset$, the Euler characteristic of the barycenter spaces has been computed for general polyhedral spaces in \cite{kk}:
\begin{equation}\label{normal}
\chi {\mathcal B}_k(X) = 1 - {k-\chi\choose k} = 1 - {1\over k!}(1-\chi)(2-\chi)\cdots (k-\chi ).
\end{equation}
It turns out that this formula is still valid for disconnected spaces
(\S\ref{disconnected}) and even if we replace $\chi$ by $\chi_c$ everywhere in the formula (Remark \ref{chicbar}).

The main contribution of this paper is to compute $\chi_c(\brho)$ for a general family of spaces $X$, not necessarily connected, and from there deduce the topological Euler characteristic $\chi$ for many cases of interest (Corollary \ref{thisisit}).

We define a \textit{basic space} to be a connected space which is either a finite CW complex or a locally closed subspace in a locally compact (Hausdorff) space. We recall that $X\subset Z$ is \textit{locally closed} (or ``LC") if it is open in its closure (see \S\ref{compact}). A typical example we consider in this paper is when $X$ is the interior of a manifold with boundary. Denote by $\mathcal P(\{1,\ldots, r\})$ (resp. $\mathcal P^*(\{1,\ldots, r\})$)
the power set of all subsets of $\{1,\ldots, r\}$ (resp. those excluding the empty set).

\bth\label{main} Let $X$ be a finite union of basic spaces and write $\chi_c=\chi_c(X)$ the Euler characteristic with compact supports of $X$. Let $p_i\in X$, $1\leq i\leq r$ be the singular points with weights $w_i>0$, and let $\rho> 0$. Then
\begin{eqnarray*}
\chi_c (\brho) &=&
1 -  {\lfloor \rho\rfloor-\chi_c +r\choose \lfloor \rho\rfloor}
- \sum_{\{i_1,\ldots, i_k\}\atop
\in{\mathcal P}^*(\{1,2,\ldots, r\})}(-1)^k{\lfloor \rho-w_{i_1}-\cdots - w_{i_k}\rfloor -\chi_c +r\choose \lfloor \rho-w_{i_1}-\cdots - w_{i_k}\rfloor}
\end{eqnarray*}
with the understanding that binomial coefficients where
$\lfloor \rho-w_{i_1}-\cdots - w_{i_k}\rfloor  < 0$ are set to zero.
\end{theorem}

Binomial coefficients can be computed in the case of negative integer entries and they have integral values (Remark \ref{negcoeffs}).
From Theorem \ref{main}, we can deduce the topological Euler characteristic $\chi \brho$ as a function of $\chi:=\chi (X)$ in the following two relevant cases.

\bco\label{thisisit} Assume $w_i\leq 1$ for all $i$. If $X$ is compact or if $X$ is the interior of an even dimensional manifold with boundary (or a union of those),
then $\chi (\brho)=\chi_c(\brho)$. In other words, the topological Euler characteristic of $\brho$ is given by the formula in Theorem \ref{main} after replacing $\chi_c$ by $\chi$ everywhere in the formula.
\eco

\begin{proof} When $X$ is compact  and the $w_i\leq 1$ for all $i$,
 $\brho$ is compact and the claim is immediate since $\chi_c$ and $\chi$ agree on compact spaces. When $X$ is the interior of a manifold with boundary $\overline{X}$, then $X$ is locally closed in $\overline{X}$, and the formula applies. If the manifold dimension is even, the Euler characteristic of the boundary is zero (being that of an odd closed manifold) and so by definition (see \S\ref{compact})
 $\chi_c(X)=\chi (\overline{X})- \chi (\partial\overline{X})= \chi (\overline{X})=\chi_c(\overline{X})$ (by compactness of $\overline{X}$). The formula in Theorem \ref{main} gives that
 $\chi_c(\brho)=\chi_c(B_\rho^{Q_r}(\overline{X}))$. By compactness of this barycenter space, this in turn is equal to $\chi({\mathcal B}_{\rho}^{Q_r}(\overline{X}))$ so that
 $$\chi_c(\brho)= \chi({\mathcal B}_{\rho}^{Q_r}(\overline{X}))$$
But a manifold with boundary $\overline{X}$ is homotopy equivalent to its interior $\overline{X}\simeq X$ via a homotopy $H$ that is supported in a collar. Since the $p_i$'s are in $X$ and can be considered to be away from the collar (after applying a homeomorphism if necessary), the homotopy $H$ can be extended to a homotopy equivalence
 ${\mathcal B}_{\rho}^{Q_r}(\overline{X}))\simeq \brho$, and the claim follows.
\end{proof}

\bre In the formula of Theorem \ref{main}, we can regard the term $\displaystyle {\lfloor \rho\rfloor-\chi_c +r\choose \lfloor \rho\rfloor}$ as the contribution of $\emptyset\in \mathcal P(\{1,\ldots, r\})$. In other words, if for $I=\{i_1,\ldots, i_k\}\subset \{1,\ldots, r\}$, we set $w_I=\sum_{j=1}^kw_{i_j}$, with
the convention that when $I=\emptyset$, $w_I=0$ and the cardinality $|I|=0$; then our formula takes the more succinct form
\begin{eqnarray*}
\chi_c (\brho) &=&
1 - \sum_{I\in{\mathcal P}(\{1,2,\ldots, r\})}(-1)^{|I|}{\lfloor \rho-w_{I}\rfloor -\chi_c +r\choose \lfloor \rho-w_{I}\rfloor}
\end{eqnarray*}
\ere

\bre\label{chicbar}
Interestingly and when there are no critical points, we see that
$\displaystyle\chi_c {\mathcal B}_k(X) = 1 - {k-\chi_c\choose k}$,
which means that the formula computing $\chi_c({\mathcal B}_k(X))$ is similar to \eqref{normal}.
\ere

\bex\label{single} When $r=1$ the sum is over $\emptyset$ and $\{1\}$, and we obtain
\begin{equation}\label{formulaq1}
\displaystyle\chi_c ({\mathcal B}_{\rho}^{Q_1}(X)) = 1 -
{\lfloor \rho\rfloor-\chi_c +1\choose \lfloor \rho\rfloor}+ {\lfloor \rho-w_{1}\rfloor -\chi_c +1\choose \lfloor \rho-w_{1}\rfloor}
\end{equation}
We can check this formula against special cases.
Write $\rho = \lfloor \rho\rfloor+\epsilon$, $0\leq\epsilon<1$.
$${\mathcal B}_{\rho}^{Q_1}(X) = \begin{cases} {\mathcal B}_{\lfloor \rho\rfloor}(X),& \ \hbox{if}\ \ \epsilon<w_1<1\ \ (\hbox{i.e.}\ \lfloor \rho-w_1\rfloor =\lfloor \rho\rfloor-1),\\
\hbox{contractible},& \ \hbox{if}\ \ w_1\leq\epsilon\ \ (\hbox{i.e.}\ \lfloor \rho-w_1\rfloor =\lfloor \rho\rfloor).\end{cases}$$
This is consistent with the Euler characteristic computation since when $\epsilon<w_1<1$,
$$\chi ({\mathcal B}_{\rho}^{Q_1}(X)) = 1 -
{\lfloor \rho\rfloor-\chi_c +1\choose \lfloor \rho\rfloor}+ {\lfloor \rho\rfloor-\chi_c\choose \lfloor \rho\rfloor-1}
= 1- {\lfloor \rho\rfloor-\chi_c\choose \lfloor \rho\rfloor}$$
and this recovers the formula in Remark \ref{chicbar}.
When $w_1\leq\epsilon$ however,
$\lfloor\rho-w_1\rfloor = \lfloor \rho\rfloor$ so that in \eqref{formulaq1},
$\chi_c{\mathcal B}_{\rho}^{Q_1}(X) =1$ always. Note that the weighted
barycenter space is contractible, but this is in general not enough to justify that $\chi_c=1$.
This peculiarity is discussed further in Example \ref{interesting}.
\eex

\bex\label{double} When $r=2$ the sum is over $\emptyset$, $\{1\}$, $\{2\}$ and $\{1,2\}$, so that
$$\chi_c( {\mathcal B}_{\rho}^{Q_2}(X)) = 1 -
 {\lfloor \rho\rfloor-\chi_c +2\choose \lfloor \rho\rfloor}
 + {\lfloor \rho-w_{1}\rfloor -\chi_c +2\choose \lfloor \rho-w_{1}\rfloor}
 + {\lfloor \rho-w_{2}\rfloor -\chi_c +2\choose \lfloor \rho-w_{2}\rfloor}
 - {\lfloor \rho-w_{1}-w_2\rfloor -\chi_c +2\choose \lfloor \rho-w_{1}-w_2\rfloor}$$
Here too, the various homotopy types for ${\mathcal B}_{\rho}^{Q_2}(X)$ can be described (see
section \ref{degree2}).
\eex

\bre To help check the validity of the formula in Theorem \ref{main}, there are two fundamental properties that must be satisfied:
\begin{itemize}
\item Under the condition $w_k\leq \rho < w_{k+1}\leq\cdots\leq w_r$ we must have
$$\chi_c (\brho) = \chi_c ({\mathcal B}_{\rho}^{Q_{k}}(X-Q_{r-k}))$$
This identity is already true at the level of spaces; i.e.
$\brho = {\mathcal B}_{\rho}^{Q_{k}}(X-Q_{r-k})$ under the stated condition.
In particular, if $\rho<w_i, \forall i$, then
$$\chi_c \brho = 1 - {\lfloor \rho\rfloor-\chi_c +r\choose \lfloor \rho\rfloor}=\chi_c ({\mathcal B}_{\lfloor \rho\rfloor}(X-Q_r))$$
The last equality follows from the fact that under the stated conditions,
$\chi_c(X\setminus Q_r) = \chi_c(X)-r$.
\item  The second fundamental property is that if $w_{i_1}=\cdots =w_{i_k}=1$, then the
points $p_{i_1},\ldots, p_{i_k}$ are not singular anymore and
$\brho = {\mathcal B}_{\rho}^{Q_{r-k}}(X)$,  so that
$$\chi_c (\brho) = \chi_c ({\mathcal B}_{\rho}^{Q_{r-k}}(X))$$
This is also verified by our formula.
\end{itemize}
\ere

The formula in Theorem \ref{main} is intimately related to the Chen-Lin degree
$d_\rho$ \cite{cl} as we mentioned earlier.

\bco\label{main1} Let $X$ and $Q_r$ as in Theorem \ref{main}.
Consider the Chen-Lin generating series
$$g(x) = (1+x+x^2+\cdots )^{-\chi_c +r}\prod_{j=1}^r(1-x^{w_j})$$
and write it in powers of $x$ as in
$$g(x) = 1 + b_1x^{n_1} + b_2x^{n_2}+\cdots + b_kx^{n_k}+\cdots $$
where $1\leq n_1<n_2<\ldots $. Suppose
$n_k\leq \rho < n_{k+1}$. Then
$$
\chi_c(\brho) = - \sum_{j=1}^kb_j = 1-d_\rho
$$
\eco

\begin{proof} One can give the proof right away and it is combinatorial.  We have the expression for $m> 0$
\begin{equation}\label{formula1}(1+x+x^2+x^3+\cdots )^{m} =
1+{m\choose 1}x+{m+1\choose 2}x^2+\cdots=
\sum_{n=0} {m+ n-1\choose n}x^n
\end{equation}
based on the  identity
\begin{equation}\label{combinatorics}
1+{m\choose 1}+{m+1\choose 2}+\cdots +{m+n-1\choose n}={m+n\choose n}
\end{equation}
Remark \ref{negcoeffs} explains why both formulas above are valid for all integers $m$.
In fact $(1+x+x^2+x^3+\cdots )^{m} = (1-x)^{-m}$ if $m$ is negative.
So starting with \eqref{formula1}, we can write
$$(1+x+x^2+x^3+\cdots )^{-\chi_c+r} =
\sum_{n=0} {-\chi_c +r+n-1\choose n}x^n$$
Multiplying this by $\prod (1-x^{w_i})$ we get the series
$g(x) = 1 + b_1x^{n_1} + b_2x^{n_2}+\cdots + b_kx^{n_k}+\cdots $.

Let $i_1\leq i_2\leq\cdots\leq i_k$ be a sequence such that
$w_{i_1}+\cdots + {w_{i_k}}\leq\rho$, and let
$i$ be the smallest integer such that
$i+w_{i_1}+\cdots + {w_{i_k}}>\rho$ (that is $i-1 = \lfloor \rho -w_{i_1}-\dots -w_{i_k}\rfloor $).
Any such sequence contributes to the coefficients of $g(x)$ via the terms with exponents
$$x^{w_{i_1}+\cdots +w_{i_k}}, x^{1+w_{i_1}+\cdots +w_{i_k}},\ldots, x^{i-1+w_{i_1}+\cdots +w_{i_k}}$$
Here the term $(-1)^kx^{w_{i_1}+\cdots +w_{i_k}}$ comes evidently from the product $\prod (1-x^{w_i})$  and the factor $x^j$ in $x^{j+w_{i_1}+\cdots +w_{i_k}}$ comes from $(1+x+x^2+x^3+\cdots )^{-\chi_c+r}$,
with the coefficient ${-\chi_c+r+j-1\choose j}$.
The total contribution from these exponents to the sum $\sum_{j=1}^k b_j$
is therefore the sum of their coefficients. To recap,
the sequence  $i_1\leq i_2\leq\cdots\leq i_k$ with $i-1=\lfloor \rho -w_{i_1}-\dots -w_{i_k}\rfloor $
contributes to $\sum_{j=1}^kb_j$ the term
$$(-1)^k\left[ 1+{-\chi_c+r\choose 1} + {-\chi_c+r+1\choose 2}+\cdots + {-\chi_c+r+i-1-1\choose i-1}\right]$$
which is equal according to \eqref{combinatorics} to
$$(-1)^k {-\chi_c+r +i-1\choose i-1}
%= (-1)^s\left[ {-\chi_c+r +i-1\choose -\chi_c+r}\right]
= (-1)^k {\lfloor \rho -w_{i_1}-\dots -w_{i_k}\rfloor -\chi_c+r \choose \lfloor \rho -w_{i_1}-\dots -w_{i_k}\rfloor}
$$
Adding these over all sequences $i_1\leq i_2\leq\cdots\leq i_k$ with the property that
$w_{i_1}+\cdots + {w_{i_k}}\leq\rho$ gives the desired identity $d_\rho=1+\sum_{j=1}^kb_j = 1-\chi_c$.
\end{proof}

\vskip 5pt
\noindent{\bf Notation}: Throughout this note we make the assumption that $\mathcal B^\emptyset = \mathcal B$,
${\mathcal B}_0(X)=\emptyset$ and that $X*\emptyset = X$.

\vskip 5pt
\noindent{\sc Acknowledgement}: This paper answers a question of A. Malchiodi who conjectured Corollary \ref{main1}. We are grateful to him for sharing his question and insight.

%%%%%%%%%%%%%%%%%%%%%%%%%%%%%%%%%%%%%%%%%%%%%%%%%%%%%%%%%%%%%%%%%%%

\section{Conical Subspaces}\label{conic}

We will make heavy use of the following construction discussed in \cite{ahmedou, carlotto}.
For $A$ closed in $X$, define
$${\mathcal B}_n(X,A) = {\mathcal B}_{n}(X)\cup \left\{\sum t_ix_i\in {\mathcal B}_{n+1}(X)\ |\ x_i\in A\ \hbox{for some $i$}\right\}$$
This is the space consisting of all configurations with at most $n$ points in $X\setminus A$ but possibly longer configurations if one of the points is in $A$.
The following is clear
\begin{itemize}
\item ${\mathcal B}_0(X,A)=A$.
\item ${\mathcal B}_n(X,A_1)\cup\cdots\cup {\mathcal B}_n(X,A_k)={\mathcal B}_n(X,\bigcup A_i)$
\end{itemize}

We can extend this definition as follows.

\bde For pairs of spaces $(X,A_i)$, define
\begin{eqnarray*}
{\mathcal B}_{n}(X,A_1,A_2,\ldots, A_k) = &\{&t_1x_1+\cdots +t_{n}x_{n} + s_1a_1 + \cdots
+ s_ka_k\ \in {\mathcal B}_{n+k}(X)\ | a_j\in \bigcup A_i\\
&& \ \sum t_i+\sum s_j=1, t_i,s_j\geq 0\}
\end{eqnarray*}
This consists of configurations in ${\mathcal B}_{n+k}(X)$ having at most $n$ points in $X-\bigcup A_i$. Again ${\mathcal B}_0(X,A_1,\ldots, A_k) = {\mathcal B}_k(\bigcup A_i)$.
These spaces are closed if the $(X,A_i)$ are closed pairs.
 By abuse of notation we write $A_i=p_i$ if $A_i=\{p_i\}$ is a singleton.
Our definition coincides with that in (\cite{cm},\S3) who adopt instead the notation $X_{i_1,\ldots, i_k}^{n,k}$ for our conical spaces $B_n(X,p_{i_1},\ldots, p_{i_k})$.
\ede

\ble\label{contractible} The conical subspaces ${\mathcal B}_{n}(X,p_1,p_2,\ldots, p_k)$
are always contractible as soon as $k\geq 1$.
\ele

\begin{proof}
We have an inclusion ${\mathcal B}_{n}(X,p_1,p_2,\ldots, p_k)\subset {\mathcal B}_{n+k-1}(X,p_1)$,
and the deformation retraction of ${\mathcal B}_{n+k-1}(X,p_1)$ to $p_1$ restricts to a deformation retraction of ${\mathcal B}_{n}(X,p_1,p_2,\ldots, p_k)$.
\end{proof}

%\item ${\mathcal B}_{n+k}(X, A_1,\ldots, A_k)\subset {\mathcal B}_{n+k}(X, A_1,\ldots,\hat{A}_i,\ldots, A_k)$, where $\hat{A}_i$ means the $i$-th entry is suppressed.
%\item If $A,B$ are disjoint, then
%${\mathcal B}_{n+1}(X,A)\cap {\mathcal B}_{n+1}(X,B)={\mathcal B}_n(X)\cup {\mathcal B}_{n+1}(X,A,B)$.
%\end{itemize}

The weighted barycenter spaces $\brho$ only depend on the homeomorphism type of $X$.
There is a major difference between the cases when the weights $w_i$ are smaller or bigger than one. In the former case, the spaces behave like ``quotients", while in the latter case they behave like ``complements". For instance, consider one singular point $p$ with weight $w$, and let $X$ be the unit disk $D$
in $\bbr^n$. If $w>1$, then
${\mathcal B}_1^{Q_1}(D) = D -\{p\}$, and the homotopy type depends on the dimension of $D$.

When the weights $w_i$ are $< 1$, it is possible to describe $\brho$ as a colimit of a diagram
of spaces of the form ${\mathcal B}_\ast(X)$ or ${\mathcal B}_\ast(X,p_{i_1},\ldots, p_{i_k})$.

\ble\label{colimit} Let $Q=\{p_1,\ldots, p_r\}$, $\Omega=\{1,\ldots, r\}$, $w(p_i)=w_i$
and suppose $0<w_i<1$ for all $i\in\Omega$. There is a colimit decomposition
\begin{eqnarray*}
\brho
= \bigcup_{\{i_1,\ldots, i_k\}\subset\Omega}{\mathcal B}_{\lfloor \rho - w_{i_1}-\cdots - w_{i_k}\rfloor }(X,p_{i_1},\ldots , p_{i_k})
\end{eqnarray*}
\ele

We can refine this union by only considering the maximal conic subspaces making up the union. As indicated by (\cite{carlotto}, Definition 2.1),
${\mathcal B}_n(X,p_{i_1},\ldots, p_{i_r})$ includes in
${\mathcal B}_m(X,p_{j_1},\ldots, p_{j_s})$ if $n\leq m$ and $\{i_1,\ldots, i_r\}$ splits
into a subset in $\{j_1,\ldots, j_s\}$ and another subset of cardinality $\leq m-n$.
Here ${\mathcal B}_n(X,p_{i_1},\ldots, p_{i_k})$ is maximal in
$\brho$ if it is not contained in a larger conic subspace.

\bex Let $\rho = 4.5$. Suppose we have $3$ singular points with weights
$w_1 = 0.3, w_2=0.4, w_3=0.6$, then
$${\mathcal B}^{Q_3}_\rho(X) = {\mathcal B}_4(X,p_1)\cup {\mathcal B}_4(X,p_2)\cup {\mathcal B}_3(X,p_1,p_2,p_3)$$
All subspaces in the union are maximal subspaces.
\eex

\bex\label{exr=1} When $r=1$, with a single singular point $p$ of weight $0<w\leq 1$, then
$${\mathcal B}^{Q_1}_\rho(X) = {\mathcal B}_{\lfloor \rho\rfloor}(X)\cup {\mathcal B}_{\lfloor \rho-w\rfloor }(X,p)$$
One can see immediately that
${\mathcal B}^{Q_1}_\rho(X) = {\mathcal B}_{\lfloor \rho\rfloor}(X)$ if $\lfloor \rho-w\rfloor <\lfloor \rho\rfloor$ since in that case $B_{\lfloor\rho-w\rfloor}(X,p)\subset B_{\lfloor\rho\rfloor}(X)$, and that
${\mathcal B}^{Q_1}_\rho(X) = {\mathcal B}_{\lfloor \rho\rfloor}(X,p)$ is contractible if $\lfloor \rho-w\rfloor =\lfloor \rho\rfloor$ (Example \ref{single}).
\eex

\bex Suppose $r=2$, $p_1,p_2$ having weights $w_1,w_2$.
In the case $w_1,w_2\leq\rho -\lfloor \rho\rfloor=\epsilon$, $w_1+w_2> \epsilon$, we can have a configuration of length $\lfloor \rho\rfloor+1$ provided
the configuration contains $p_1$ or $p_2$, but no configuration can be of length $\lfloor \rho\rfloor+2$. This means
that
$$B^{Q_2}_\rho(X) = {\mathcal B}_{\lfloor \rho\rfloor}(X,p_1)\cup {\mathcal B}_{\lfloor \rho\rfloor}(X,p_2)$$
Other descriptions occur depending on the choices of $w_1,w_2$
(see section \ref{degree2}).
Notice that in the present case
$$ {\mathcal B}_{\lfloor \rho\rfloor}(X,p_1)\cap {\mathcal B}_{\lfloor \rho\rfloor}(X,p_2) = {\mathcal B}_{\lfloor \rho\rfloor}(X)\cup {\mathcal B}_{\lfloor \rho\rfloor-1}(X,p_1,p_2)$$
It is always true that the intersection of conical subspaces is again a union of conical subspaces. This fact is crucial when it comes to determining some quotients and some homotopy types.
\eex

We next list the main properties for the conic spaces needed in our
computation of the Euler characteristic.
All spaces below are connected, hausdorff and locally compact.

\ble\label{properties} Assume $n\geq 1$. The following hold\\
$\bullet$ (i)  ${\mathcal B}_n(X,A)$ is contractible if $A$ is contractible. In particular
${\mathcal B}_n(X,p)$ is contractible.\\
$\bullet$  (ii) \cite{ahmedou} Let $A$ be a closed subspace of $X$. Then
$\displaystyle {{\mathcal B}_n(X)\over {\mathcal B}_{n-1}(X,A)} \simeq {\mathcal B}_n(X/A)$\\
$\bullet$ (iii) If $X$ is contractible, then
$\displaystyle {\mathcal B}_n(X/A)\simeq \Sigma {\mathcal B}_{n-1}(X,A)$ where $\Sigma$ means suspension.
\ele

\begin{proof} (i) is Lemma \ref{contractible}.
For the second claim, notice that
 $\displaystyle {{\mathcal B}_n(X)\over {\mathcal B}_{n-1}(X,A)} = {{\mathcal B}_n(X/A)\over
{\mathcal B}_{n-1}(X/A, p)}$ where $p$ is the preferred basepoint in $X/A$.
Since ${\mathcal B}_{n-1}(X/A, p)$ is contractible we obtain the homotopy equivalence in claim (ii).
The simplest example here is when $n=1$, $B_0(X,A)=A$, $B_1(X)=X$ and the quotient is $X/A$.
Claim (iii) is a direct consequence of (ii) and the fact that ${\mathcal B}_n(X)$ is also contractible.
\end{proof}

\bex\label{xu}
When $X=D$ is a closed $m$-dimensional ball with boundary $\partial D=S^{m-1}$,
${\mathcal B}_n(X)$ is contractible since $D$ is contractible, so
lemma \ref{properties} (iii) implies that $ {\mathcal B}_{n}(S^m)\simeq \Sigma ({\mathcal B}_{n-1}(D,
S^{m-1}))$. Note that $D/\partial D=S^m=S^m/*$, but ${\mathcal B}_n(D,\partial D)$ behaves quite differently than
${\mathcal B}_n(S^m,p)$ which is contractible.
\eex

%%%%%%%%%%%%%%%%%%%%%%%%%%%%%%%%%%%%%%%%%%%%%%%%%%%%%%%%%%%%%%%%%%%

\section{The Compactly Supported Euler Characteristic}\label{compact}

The compactly supported Euler characteristic $\chi_c$, sometimes called the ``combinatorial" Euler characteristic, is defined for locally compact spaces and has the property that for any
disjoint decomposition of $X=\coprod X_i$ where each $X_i$ is a locally closed subspace of $X$,
$$\chi_c(X) = \sum\chi_c(X_i).$$
As is common, we reserve the word ``stratification" $\{X_i\}$ for $X$ if $X$ is a disjoint union of the $X_i$'s and all $X_i$'s are locally closed in $X$.

\bre Being locally closed in a topological space $X$ has various equivalent definitions (see \cite{gr}, Proposition 1). Here's conveniently the list of equivalences: $A$ is locally closed in $X$ $\Longleftrightarrow$ $A= U\cap cl A$ for some open set $U$
$\Longleftrightarrow$ $cl A\setminus A $ is closed $\Longleftrightarrow$ $A\cup (X\setminus cl A)$ is open. It seems more convenient to think in terms of the third equivalence: $A$ is LC in $X$ if $cl A\setminus A$ is closed.
\ere

The above additivity formula for $\chi_c$ makes it a very computable characteristic. Its drawback is that it is not an invariant of homotopy type. For example when $D$ is the open unit disk, then
$\chi (D)=1$ if $D$ is even dimensional, and $\chi (D)=-1$ if $D$ has odd dimension.
In particular, if $X$ is contractible, $\chi_c(X)$ is not necessarily $1$.

As expected, $\chi_c(X)=\chi (X)$ if $X$ is compact.

We now explain how to compute $\chi_c$.
Borel-Moore homology can be computed as follows (our main reference is
\cite{ginzburg}, \S2.6, also \cite{goresky}). Let $X$ be a locally compact hausdorff space with compactification $\overline{X}$ of $X$ such that $\overline{X}\setminus X$ is closed in $\overline{X}$. Then (with field coefficients)
$$H_*^{BM}(X) \cong H_*(\overline{X},\overline{X}\setminus X)$$
and this is valid for disconnected spaces as well. Moreover, it is known that
$H_*^{BM}(X)\cong H_c^*(X)$ the cohomology with compact support with field coefficients.
Let $\chi_c(X)$ be the alternating sum of the ranks of these groups.
For a pair $(\overline{X},X)$ as above, we must have
\begin{equation}\label{mainchic}
\chi_c(X) =  \chi (\overline{X})
-\chi (\overline{X}\setminus X) = \chi \left({\overline{X}\over\overline{X}\setminus X}\right)-1
\end{equation}
The $1$ is subtracted to take into account the loss of one factor in $H_0$. This is the formula we use to compute $\chi_c$ throughout the paper.

\bde\label{hadhi} As in \cite{kt}, we say that $X$ admits a ``BM-compactification" if $X$ is locally compact hausdorff with compactification $\overline{X}$ such that
$\overline{X}\setminus X$ is closed in $\overline{X}$.
If $X$ is compact, then
$\overline{X}=X$.
\ede

\bex \cite{kt} If $X,Y$ admit BM-compactifications, then
\begin{equation}\label{chijoin}
\chi_c(X*Y) = \begin{cases} - \chi_c(X)\chi_c(Y)&\ \hbox{$X$ and $Y$ both not compact.}\\
\chi_c (X)+\chi_c (Y)-\chi_c (X)\chi_c (Y)&
\ \hbox{either $X$ or $Y$ is compact, or both}.\\ \end{cases}
\end{equation}
For disks, the above is consistent with the fact that $D^n*D^m\cong D^{n+m+1}$ and
 $\chi_c(D^n)=(-1)^n$. Also and as a consequence of this formula, we
 find that
\begin{eqnarray}
\chi_c (\Sigma^kX) &=& 1+(-1)^k(\chi_c (X)-1)\label{susp}
\end{eqnarray}
where $\Sigma^kX:=S^k*X$ is the suspension iterated $k$ times of $X$.
When $k$ is even $\chi_c (\Sigma^kX)=\chi_cX$ and when $k$ is odd
$\chi_c (\Sigma^kX) =2-\chi_c (X)$.
\eex

\iffalse
The next result is of course straigthforward using the additivity of $\chi_c$ but the method of calculation is instructive and illustrative of our methods.

\ble\label{illustrative} Given a BM compactification $(\overline{X},\overline{X}\setminus X)$, $Q_r$ a set of $r$-points in $X$, then $\chi_c(X-Q_r)=\chi_c(X) -r$.
\ele

\begin{proof}
We want to see how to obtain this formula using the definition in \eqref{mainchic}.
We use the crucial fact, valid in simplicial complexes,
that any point $p\in X$ has a open neighborhood $U$, which is a deformation retract, and such that $X\setminus\overline{U}\cong X\setminus\{p\}$. We say such a neighborhood is ``adapted".
Replace $X-Q_r$ by
$X-\bigcup^r U_i$ where each $U_i$ is an adapted open neighborhood around $p_i$, the neighborhoods do not intersect. A compactification of $X-\bigcup^r\overline{U}_i$ is $X-\bigcup^r U_i$.
This is a BM-compactification with difference the boundary
$\bigcup^r\partial\overline{U}_i$. Then
$$\chi_c (X-Q_r)=\chi \left( {X-\bigcup^r U_i\over \bigcup^r\partial\overline{U}_i}\right)-1=
\chi (X\vee \bigvee^{r-1} S^1)-1 = \chi (X) - r$$
\end{proof}

\fi

\bex\label{interesting} This next example is pertinent and discusses the computation of $\chi_c(B_\rho^{Q_1}(D))$, where $D$ is the \textit{open} unit disk in $\bbr^n$, $\rho = \lfloor\rho\rfloor +\epsilon$ and where the unique singular point $p_1$ at the origin has weight $0<w_1\leq\epsilon$. This is a good illustration of the importance of having the locally closed condition when computing $\chi_c$, and is also some explanation of the peculiarity discussed at the end of Remark \ref{single}. Now evidently in this case ${\mathcal B}_\rho^{Q_1}(D)={\mathcal B}_{\lfloor\rho\rfloor}(D,p_1)$ is contractible. Take $\lfloor\rho\rfloor = 1$, then ${\mathcal B}_{1}(D,p_1)$ consists of all configurations of the form
$t_1x+t_2p_1$, $t_1+t_2=1$. This is the cone on $D$ but we must identify all points of the form $t_1p_1+t_2p_1\sim p_1$, so ${\mathcal B}_{1}(D,p_1)$ is homeomorphic to the \textit{reduced cone
$cD$} obtained from the cone by collapsing the segment $[v,0]$ to $0$, where $v$ is the vertex of the cone, and $0$ the origin of $D$. This space is stratified as follows. We will specify the dimension by writing $D^n$ for the $n$-dimensional open disk. Its cone with vertex $v$ can be stratified as $ D^{n+1}\sqcup \{v\}\sqcup D^n$. The reduced cone $cD$ has then a locally closed stratification of the form (up to homeomorphism)
$$cD = (D^{n+1}\setminus L)\sqcup D^n$$
where $L$ is a diameter in the open disk $D^{n+1}$ (note already that $v$ has been identified with $0\in D^n$ in this decomposition). A BM-compactification of $D^{n+1}\setminus L$ is the closed disk $\overline{D}^{n+1}$
and $\overline{D}^{n+1}\setminus (D^{n+1}\setminus L)$ is the sphere $S^n$ together with a diameter attached, so it is homotopy equivalent to $S^n\vee S^1$.
We can then write
\begin{eqnarray*}
\chi_c({\mathcal B}_{1}(D,p_1))=\chi_c(cD)&=&
\chi_c({D}^{n+1}\setminus L)+\chi_c(D^n)\\
&=&\chi (\overline{D}^{n+1})-\chi (S^n\vee S^1) + \chi (\overline{D}^n)-\chi (S^{n-1})\\
&=&1- (\chi (S^n)-1) + 1 -\chi (S^{n-1}) = 1
\end{eqnarray*}
so this is always $1$ as predicted by Theorem \ref{main}. The way we compute this in
\S\ref{mainkey} is by using a slightly different stratification which is better adapted
to the general situation.
\eex

%%%%%%%%%%%%%%%%%%%%%%%%%%%%%%%%%%%%%%%%%%%%%%%%%%%%%%%%%%%%%%%%%%%%%%%%%%%%%%%%

\subsection{Topology}\label{adjunct}

The topology of $\brho$ can be subtle, especially when dealing with pushouts and compactifications.  We make below a few relevant observations related to these issues.

If $A\subset X$ is a closed subspace, then
${\mathcal B}_n(A)$ is closed in ${\mathcal B}_n(X)$. This is not generally true if we replace closed by open. Take $A=(0,1)$ and $X=[0,1]$. The sequence $((1-{1\over n})x_1 + {1\over n}x_2)$, with $x_1={1\over 2}$ and $x_2=1$, is in
${\mathcal B}_2([0,1])\backslash {\mathcal B}_2((0,1))$ and converges to $x_1\in (0,1)={\mathcal B}_1(0,1)\subset {\mathcal B}_2(0,1)$, which means that
${\mathcal B}_2(A)$ cannot be open in ${\mathcal B}_2(X)$. That the statement is true for closed sets and not open sets is a consequence of the fact that generally
${\mathcal B}_n(X)\setminus {\mathcal B}_n(A)$ is not ${\mathcal B}_n(X\setminus A)$. In fact and
\textit{as sets}
$${\mathcal B}_n(X) = {\mathcal B}_n(X\setminus A)\cup {\mathcal B}_{n-1}(X,A)$$
This union is not a pushout or an adjunction space construction (\cite{brown}, section 4.5).
Indeed the intersection of both factors in the union is
${\mathcal B}_{n-1}(X\setminus A)$. The union doesn't have the quotient topology of the corresponding adjunction space
since there is part of the boundary of ${\mathcal B}_n(X\setminus A)$ intersecting ${\mathcal B}_{n-1}(X,A)$ and yet not being in ${\mathcal B}_{n-1}(X\setminus A)$.

\begin{notation}\label{notation} \rm (and terminology) Recall that $X$ is basic if it is connected, either compact or locally closed in a locally compact space, with closure $\overline{X}$. If $X$ is compact, then we write $\overline{X}=X$, and if not, then $\overline{X}$ is the closure of $X$, and
$(\overline{X},X)$ is a BM-compactification. In this latter case, we
write conveniently $\partial\overline{X}:=\overline{X}\setminus X$.
\end{notation}

\ble\label{keyprime} Let $Q_r=\{p_1,\ldots, p_r\}\subset X$,
and suppose $X$ is basic. Then ${\mathcal B}_n(X-Q_r)-{\mathcal B}_{n-1}(X-Q_r)$ is locally closed in ${\mathcal B}_n(\overline{X})$, and this pair is a BM-compactification.
\ele

\begin{proof}
The point is that a configuration approaching the closure of
${\mathcal B}_n(X-Q_r)-{\mathcal B}_{n-1}(X-Q_r)$ in ${\mathcal B}_n(\overline{X})$ must have a point approaching $Q_r\cup\partial\overline{X}$,
or two points of the configuration approach each other.
In other words and more precisely, we have
\begin{eqnarray*}
{\mathcal B}_n(\overline{X})\setminus ({\mathcal B}_n(X-Q_r)-{\mathcal B}_{n-1}(X-Q_r))
&=& ({\mathcal B}_n(\overline{X})\setminus {\mathcal B}_n(X-Q_r))\cup {\mathcal B}_{n-1}(\overline{X}-Q_r)\label{says} \\
&=& {\mathcal B}_{n-1}(\overline{X},Q_r\cup\partial\overline{X})
\end{eqnarray*}
This is closed in ${\mathcal B}_n(\overline{X})$ which is compact and is the closure of
${\mathcal B}_n(X-Q_r)-{\mathcal B}_{n-1}(X-Q_r)$.
\end{proof}

\iffalse
\ble\label{key} Suppose $A\subset X$ open, $X$ compact. Then
${\mathcal B}_n(X-\bar A)-{\mathcal B}_{n-1}(X-\bar A)$ is locally closed in ${\mathcal B}_n(X-A)$.
\ele

\begin{proof}
For $n=0$, ${\mathcal B}_0(X)$ was declared empty and ${\mathcal B}_0(X,A)=A$. The claim is clear. Suppose
 $n\geq 1$ and write
\begin{eqnarray}
{\mathcal B}_n(X-A)\setminus ({\mathcal B}_n(X-\bar A)-{\mathcal B}_{n-1}(X-\bar{A}))
&=& ({\mathcal B}_n(X-A)\setminus {\mathcal B}_n(X-\bar A))\cup {\mathcal B}_{n-1}(X-\bar{A})\label{says} \\
&=& {\mathcal B}_{n-1}(X- A,\partial\overline{A})\nonumber
\end{eqnarray}
and this is closed in ${\mathcal B}_n(X-A)$.
Since $X\setminus A$ is closed in $X$, it is thus compact. Thus ${\mathcal B}_n(X-A)$ is compact
and \eqref{says} says that it is
a BM-compactification of ${\mathcal B}_n(X-\bar A)-{\mathcal B}_{n-1}(X-\bar A)$.
\end{proof}
\fi

The lemma above is what enables us to compute
$\chi_c$ for such complements in section \ref{relative} next.

Another remark pertaining to topology is to take $A,B$ disjoint in $X$.
For example can have $B=X\setminus A$. Then
${\mathcal B}_i(A)*{\mathcal B}_j(B)$ is a subset of ${\mathcal B}_{i+j}(X)$, meaning the topology of the join coincides with the induced topology. The closure of ${\mathcal B}_i(A)*{\mathcal B}_j(B)$ in ${\mathcal B}_{i+j}(X)$ is homeomorphic to ${\mathcal B}_i(\overline{A})*{\mathcal B}_j(\overline{B})$ if and only if however the closures are disjoint $\overline{A}\cap\overline{B}=\emptyset$.

\bex\label{lcstratum} Consider the stratum
in $\brho$ consisting of all configurations having
\textit{exactly} $k$ singular points $p_1,\ldots, p_k$ appearing in the configuration (i.e. having non-zero coefficients), and no other singular points. Let's denote this stratum by
${\mathcal B}_n(X-Q_r)\overset{\circ}{*}{\mathcal B}_k\{p_1,\ldots, p_k\}$.
A configuration in  ${\mathcal B}_n(X-Q_r)\overset{\circ}{*}{\mathcal B}_k\{p_1,\ldots, p_k\}$ is then of the form $\sum_{i=0}^n t_ix_i+ s_1p_1+\cdots +s_kp_k, s_i\neq 0$, $x_i\not\in Q_r$. This stratum as a subspace
of $\brho$ can be written as
\begin{equation} \label{stratum}
{\mathcal B}_n(X-Q_r)\overset{\circ}{*}{\mathcal B}_k\{p_1,\ldots, p_k\}:=
{\mathcal B}_n(X-Q_r)*{\mathcal B}_{k}\{p_1,\ldots, p_k\} -
{\mathcal B}_n(X-Q_r)*{\mathcal B}_{k-1}\{p_1,\ldots, p_k\}
\end{equation}
Its BM-compactification can therefore be deduced from Lemma \ref{keyprime} and \eqref{chijoin} (see Proposition \ref{secondchic}).
\eex

%%%%%%%%%%%%%%%%%%%%%%%%%%%%%%%%%%%%%%%%%%%%%%%%%%%%%%%%%
\section{Proof of Theorem \ref{main}}\label{relative}

This proof is broken in several steps. We
recall throughout that ${\mathcal B}_0(X)=\emptyset$ and that $X*\emptyset = X$. Notation
and terminology are as in Definition \ref{hadhi} and Notation \ref{notation}.

\ble\label{firstchic} Let $X$ be a basic set. Then
$\chi_c ({\mathcal B}_n(X))=1-{n-\chi_c\choose n}$.
\ele

\begin{proof}
The proof is carried out for $(\overline{X},X)$ pair, the case when $X$ compact is known.
The BM-compactification of ${\mathcal B}_n(X)$ is certainly not ${\mathcal B}_n(\overline{X})$ as indicated in Lemma \ref{keyprime}. We therefore need to stratify ${\mathcal B}_n(X)$ as follows. Set
$$X_n := {\mathcal B}_n(X)-{\mathcal B}_{n-1}(X)$$
so that in light of the aforementioned lemma, $\overline{X}_n = {\mathcal B}_n(\overline{X})$ and
$\overline{X}_n\setminus X_n = {\mathcal B}_{n-1}(\overline{X},\partial\overline{X})$,
where this last term we recall consists of all
configurations $\sum t_ix_i\in {\mathcal B}_n(\overline{X})$ with $x_i\in \partial\overline{X}$ for some $i$ or $t_i=0$ for some $i$. We have
\begin{eqnarray*}
{\overline{X}_n\over \overline{X}_n\setminus X_n} =
{{\mathcal B}_n(\overline{X})\over {\mathcal B}_{n-1}(\overline{X},\partial\overline{X})}
&\simeq&{{\mathcal B}_n\left({\overline{X}\over \partial\overline{X}}\right)} \ \ \ \ \ \ \ \hbox{by Lemma \ref{properties}}
\end{eqnarray*}
and we get immediately that for $n\geq 2$
\begin{equation}\label{mallafaza}
\chi_c(X_n) = \chi (\overline{X}_n/ \overline{X}_n\setminus X_n) - 1
= \chi\left({\mathcal B}_n\left({\overline{X}\over\partial\overline{X}}\right)\right) -1 =
-{n-\chi_c-1\choose  n}
\end{equation}
with the last equality obtained from \eqref{normal} and the identity \eqref{mainchic}.

To get to $\chi_c({\mathcal B}_n(X))$ we write
\begin{eqnarray*}
{\mathcal B}_n(X) &=& \left({\mathcal B}_n(X)\setminus {\mathcal B}_{n-1}(X)\right)\sqcup \cdots
\sqcup ({\mathcal B}_2(X)\setminus X) \sqcup X\nonumber \\
&\cong&X_n\sqcup X_{n-1}\sqcup\cdots\sqcup X_1
\end{eqnarray*}
which is a stratification with locally closed strata, so
we can use the additivity of $\chi_c$ to see that
\begin{eqnarray*}
\chi_c{\mathcal B}_n(X) =\sum_{1\leq i\leq n} \chi_c(X_i)
&=& -{n-\chi_c-1\choose n} - {n-\chi_c-2\choose n-1}-\cdots - {-\chi_c\choose 1}\\
&=&1-{n-\chi_c\choose n}
\ \ \ \ \hbox{by}\ \ \eqref{combinatorics}
\end{eqnarray*}
which proves the lemma.
\end{proof}

\bco Let $X$ be a basic set and $Q_r\subset X$. Then
$\displaystyle\chi_c ({\mathcal B}_n(X-Q_r))=1-{n-\chi_c +r\choose n}$.
\eco

The next two lemmas are a preparation for Proposition \ref{secondchic}.

\ble\label{identification} For $X$ basic and $Q_r\in X$, we have
$${\overline{X}\over Q_r\cup\partial\overline{X}}\simeq
\begin{cases} {\overline{X}\over\partial\overline{X}}\vee\bigvee^{r} S^1& \hbox{if}\
\partial\overline{X}\neq\emptyset\\
X\vee\bigvee^{r-1} S^1,& \hbox{if $X$ compact}.
\end{cases}$$
\ele

This a straightforward consequence of the fact that in a path-connected space, the identification of two points is up to homotopy like the one point union with a circle. The proof is skipped. The analog of this lemma for disconnected spaces is given in Lemma \ref{quotient}. The next Lemma is imported from \cite{ahmedou}.

\ble\label{ahmadoudou} (\cite{ahmedou} Lemma 5.10) Let $(X,A)$ and $(Y,B)$ be two connected CW pairs. Then $(X*Y)/(A*Y)\simeq \left(X/A\right)* Y$ and
$$(X*Y)/(X*B\cup A*Y)\simeq
\begin{cases} X/A *Y/B &,\  A\neq\emptyset, B\neq\emptyset,\\
\Sigma (X/A) \rtimes  Y &,\  A\neq\emptyset, B=\emptyset,\\
\Sigma (X\times Y)\vee S^1&,\ A=\emptyset, B=\emptyset.
\end{cases}
$$
where $X\rtimes Y:= {X\times Y\over x_0\times Y}$ is the half-smash product.
\ele

\bpr\label{secondchic} Again $X$ basic, $Q_r\subset X$ and $\chi_c:=\chi_c(X)$.
For $k\geq 1$, define  ${\mathcal B}_n(X-Q_r)\overset{\circ}{*}{\mathcal B}_k\{p_1,\ldots, p_k\}$ as in \eqref{stratum}. Then
$$\chi_c \left({\mathcal B}_n(X-Q_r)\overset{\circ}{*}{\mathcal B}_k\{p_1,\ldots, p_k\}\right)= (-1)^{k+1}{n-\chi_c +r\choose n}$$
\epr

\begin{proof}
We recall that
${\mathcal B}_kQ_k= {\mathcal B}_k\{p_1,\ldots, p_k\}$ is (i.e. homeomorphic to)
 the $k-1$ dimensional simplex $\Delta_{k-1}$ and
that ${\mathcal B}_{k-1}Q_k$ is its boundary sphere $\partial\Delta_{k-1}$.
Consider as in Example \ref{lcstratum} the subspace
\begin{eqnarray*}\label{secondstrat}
X_i &=& {\mathcal B}_i(X-Q_r)\overset{\circ}{*}{\mathcal B}_kQ_k -
{\mathcal B}_{i-1}(X-Q_r))\overset{\circ}{*}{\mathcal B}_kQ_k
\end{eqnarray*}
We clearly have
$${\mathcal B}_n(X-Q_r)\overset{\circ}{*}{\mathcal B}_kQ_k =
X_n\sqcup X_{n-1}\sqcup\cdots\sqcup X_1\sqcup X_0$$
where the last two spaces are given by
$X_1 = (X- Q_r)\overset{\circ}{*}{\mathcal B}_kQ_k - \overset{\circ}{B}_kQ_k$
and $X_0 = \overset{\circ}{B}_kQ_k = {\mathcal B}_kQ_k- {\mathcal B}_{k-1}Q_k$ (the interior of a disk of dimension $k-1$). Note that the first stratum $X_n$ is open in ${\mathcal B}_n(X-Q_r)\overset{\circ}{*}{\mathcal B}_kQ_k$.
This is a locally closed stratification and
$\chi_c({\mathcal B}_n(X-Q_r)\overset{\circ}{*}{\mathcal B}_kQ_k=\sum_{0\leq i\leq n}\chi_c(X_i)$.

The simplest case is when $n=0$ where evidently
\begin{equation}\label{firstterm}
\chi_c(X_0)=\chi_c(\overset{\circ}{\mathcal B}_kQ_k) = \chi_c(\overset{\circ}{\Delta}_{k-1})= (-1)^{k-1}
\end{equation}

For $n\geq 1$, and assuming $X$ is compact, a BM-compactification for $X_n$ is given by
$\overline{X}_n={\mathcal B}_n(X)*{\mathcal B}_k(Q_k)$.  The complement of $X_n$ in $\overline{X}_n$ is
$$({\mathcal B}_n(X)*{\mathcal B}_{k-1}(Q_k))\cup ({\mathcal B}_{n-1}(X,\{p_1,\ldots, p_r\})*{\mathcal B}_k (Q_k))$$
so that from the definition of $\chi_c$ we have
$$\chi_c(X_n) =\chi\left({{\mathcal B}_n(X)*{\mathcal B}_k(Q_k)\over
({\mathcal B}_n(X)*{\mathcal B}_{k-1}(Q_k))\cup ({\mathcal B}_{n-1}(X,\{p_1,\ldots, p_r\})*{\mathcal B}_k (Q_k))}\right)-1$$
The quotient in the first term on the right can be described precisely
(Lemma \ref{ahmadoudou}). For $k\geq 2$, this is
\begin{eqnarray}\label{several}
\chi_c(X_n) &=&\chi\left({{\mathcal B}_n(X)\over
{\mathcal B}_{n-1}(X,\{p_1,\ldots, p_r\})}*
{{\mathcal B}_kQ_k\over {\mathcal B}_{k-1}Q_k}\right)-1\ \ \ \ \ \ \ \ \ \ \ \hbox{by Lemma \ref{ahmadoudou}, first case}\\
&=&\chi \left({\mathcal B}_n(X\vee\bigvee^{r-1}S^1)*S^{k-1}\right)-1\nonumber\\
&=&\chi\left(\Sigma^k{\mathcal B}_n(X\vee\bigvee^{r-1}S^1)\right)-1\nonumber\\
&=&(-1)^k\left(\chi ({\mathcal B}_n (X\vee\bigvee^{r-1} S^1)-1\right)\ \ \ \ \ \ \ \ \ \ \ \hbox{by \eqref{susp}}\nonumber\\
&=&-(-1)^k{n-\chi+r-1\choose n}\nonumber\ \ \ \ \ \ \ \ \ \ \ \ \ \ \ \ \ \ \ \ \ \ \hbox{by \eqref{normal}}\label{secondterm}
\end{eqnarray}
Adding up \eqref{firstterm} and\eqref{secondterm}, we obtain (in case $k\geq 2$)
\begin{eqnarray*}
&&\chi_c ({\mathcal B}_n(X-Q_r)\overset{\circ}{*}{\mathcal B}_k(\{p_1,\ldots, p_k\})\\
&&= -(-1)^k\left({n-\chi+r-1\choose n}+{n-1-\chi+r-1\choose n-1}+\cdots +{1-\chi+r-1\choose 1} + 1\right)
\end{eqnarray*}
and this is $\displaystyle -(-1)^k{n-\chi+r\choose n}$ as desired.
This covers the case $k\geq 2$.

For the case $k=1$, the exact same steps as in \eqref{several} apply only
that the first step becomes
$$\chi_c(X_n)=\chi\left({B_n(X)*p\over B_n(X)\cup B_{n-1}(X,Q_r)*p }\right)-1=\chi\left(\Sigma B_n\left({X\over Q_r}\right)\right)-1$$
by the second case of Lemma \ref{ahmadoudou}. The rest is the same.

To double check the case $k=1$  using Lemma \ref{keyprime}, we can write
${\mathcal B}_n(X-Q_r)\overset{\circ}{*}p = {\mathcal B}_n(X-Q_r)*p - {\mathcal B}_{n}(X-Q_r)$,
so we have a similar stratification by locally closed subsets
\begin{eqnarray*}
X_i &=& ({\mathcal B}_i(X-Q_r)*p\setminus {\mathcal B}_i(X-Q_r)) -
({\mathcal B}_{i-1}(X-Q_r)*p\setminus {\mathcal B}_{i-1}(X-Q_r))\\
&=& ({\mathcal B}_i(X-Q_r)-{\mathcal B}_{i-1}(X-Q_r))*p - ({\mathcal B}_i(X-Q_r)-{\mathcal B}_{i-1}(X-Q_r)) -p
\end{eqnarray*}
for $i\geq 1$ and $X_0 = p$. As before
${\mathcal B}_n(X-Q_r)\overset{\circ}{*}p = X_n\sqcup X_{n-1}\sqcup\cdots\sqcup X_1\sqcup X_0$
and we can compute the $\chi_c$ of each stratum.
Since $\chi_c(Y*p) = 1$ if $Y$ has a BM-compactification,
see \eqref{chijoin}, we have that $\chi_c(Y*p-p)=0$.
This gives that
\begin{eqnarray*}
\chi_c({\mathcal B}_n(X-Q_r)\overset{\circ}{*}p) &=& -\sum_{1\leq i\leq n}\chi_c({\mathcal B}_i(X-Q_r)-{\mathcal B}_{i-1}(X-Q_r)) + \chi (p)\\
&=&1 + \sum_{1\leq i\leq n}{i-\chi+r-1\choose i} \ \  \ \ \ \ \ \hbox{by \eqref{mallafaza}}\\
&=& {n-\chi+r\choose n}
\end{eqnarray*}
which is what we wanted to prove.

In conclusion, the above calculations give the right answer in the case $X$ is compact. When $(\overline{X},X)$ is a BM compactification, we have to modify each step of the proof by incorporating a boundary $\partial\overline{X}$, keeping in mind that the singular points $p_j$ are in the interior. The modifications are indicated below, yielding the same final results. For $k\geq 2$, we use the same stratification by $X_n$, only that $\overline{X}_n$ will change.
\begin{itemize}
\item For $n\geq 1$,
$\overline{X}_n={\mathcal B}_n(\overline{X})*{\mathcal B}_k(Q_k)$, and the complement
$\overline{X}_n\setminus X_n$ is
$$({\mathcal B}_n(\overline{X})*{\mathcal B}_{k-1}(Q_k))\cup ({\mathcal B}_{n-1}(\overline{X},\partial\overline{X}\cup\{p_1,\ldots, p_r\})*{\mathcal B}_k (Q_k))$$
\item When $k\geq 2$, the quotient ${\overline{X}_n\over \overline{X}_n\setminus X_n}$ is up to homotopy ${\mathcal B}_n\left({\overline{X}\over \partial\overline{X}}\vee\bigvee^rS^1\right)*S^{k-1}$. Lemma \ref{identification}, and the same computation as in \eqref{secondterm} yields the same formula with $\chi$ replaced by $\chi_c$.
\item For $k=1$, no changes, simply apply \eqref{mallafaza} with $\chi_c$ instead of $\chi$.
\end{itemize}
This concludes the proof if $X$ is a basic set.
\end{proof}

%%%%%%%%%%%%%%%%%%%%%%%%%%%%%%%%%%%%%%%%%%%%%%%%%%%%%%%%%%%%

\subsection{Proof of the Main Theorem}\label{mainkey}

We derive the formula of Theorem \ref{main} for a basic set $X$ and singular points $p_1,\ldots, p_r$. We will prove that
\begin{eqnarray*}
\chi_c (\brho) &=&
1 -  {\lfloor \rho\rfloor-\chi_c +r\choose \lfloor \rho\rfloor}
- \sum_{\{i_1,\ldots, i_k\}\atop
\in{\mathcal P}^*(\{1,2,\ldots, r\})}(-1)^k{\lfloor \rho-w_{i_1}-\cdots - w_{i_k}\rfloor -\chi_c +r\choose \lfloor \rho-w_{i_1}-\cdots - w_{i_k}\rfloor}
\end{eqnarray*}

\begin{proof}
The key point is to write
$\brho$ as the \textit{disjoint} union of subspaces
\begin{eqnarray}\label{disjointunion}
\brho&=& {\mathcal B}_{\lfloor \rho\rfloor}(X-Q_r)\sqcup
\coprod_i {\mathcal B}_{\lfloor \rho-w_i\rfloor }(X-Q_r)\overset{\circ}{*}p_i\ \sqcup \\
&&\coprod_{\{i_1, i_2\}}{\mathcal B}_{\lfloor \rho-w_{i_1}-w_{i_2}\rfloor }(X-Q_r)\overset{\circ}{*}{\mathcal B}_2\{p_{i_1},p_{i_2}\}\sqcup\coprod_{\{i_1, i_2, i_3\}}\cdots\nonumber\\
&&\cdots \sqcup\ {\mathcal B}_{\lfloor \rho-w_{1}-\cdots - w_{r}\rfloor }(X-Q_r)\overset{\circ}{*}{\mathcal B}_r\{p_1,\ldots ,p_r\}\nonumber
\end{eqnarray}
In this notation ${\mathcal B}_i(X-Q_r)\overset{\circ}{*}Z$ is empty if $i<0$,
${\mathcal B}_0(X-Q_r)\overset{\circ}{*}Z =Z$
 and $i_n\neq i_m$ if $n\neq m$.
In this disjoint union, notice that if $p_i$ has weight $w_i\leq \rho$, then it appears in the term
${\mathcal B}_{\lfloor \rho-w_i\rfloor }(X-Q_r)\overset{\circ}{*}p_i$ (and it appears only there) so everything is accounted for once. We claim that
$$\chi_c (\brho) = \chi_c({\mathcal B}_{\lfloor \rho\rfloor}(X-Q_r))+\sum_{\{i_1,\ldots, i_k\}}\chi_c\left(
{\mathcal B}_{\lfloor \rho-w_{i_1}-\cdots - w_{i_k}\rfloor }(X-Q_r)\overset{\circ}{*}{\mathcal B}_2\{p_{i_1},\cdots, p_{i_k}\}\right)
$$
and this sum yields precisely the desired formula in light of Proposition \ref{secondchic} and lemma \ref{firstterm}. As pointed out multiple times, the factors
${\mathcal B}_{\lfloor \rho-w_{i_1}-\cdots - w_{i_k}\rfloor }(X-Q_r)\overset{\circ}{*}{\mathcal B}_2\{p_{i_1},\cdots, p_{i_k}\}$
are not locally closed in $\brho$ in general, but they are themselves stratified by strata $\{X_i\}$ which are LC in $\brho$, so $\chi_c$ is additive on \eqref{disjointunion} and the proof follows.
\end{proof}

\bre\label{negcoeffs}
We point out a computational aspect of this formula.
Binomial coefficients $\displaystyle{n\choose k}$ can be computed in the case of negative integers $n$, and
non-negative $k$. One should view $\displaystyle {X\choose k}$ as the rational function
$\displaystyle {X(X-1)\cdots (X-(k-1))\over k(k-1)\cdots 1}$, so
substituting $n$ (any integer) for $X$ gives that
for $0\leq n< k$, $\displaystyle {n\choose k} = 0$, while for $n< 0$ (eg. \cite{kronenburg})
$${n\choose k} = (-1)^k{-n+k-1\choose k}\ \ ,\ \ k > 0, n< 0$$
We set ${n\choose 0}=1$ for all $n\in\mathbb Z$. Obviously we can check that with this definition, the identity in \eqref{combinatorics} remains valid. This is because it is valid at the level of rational functions (replacing $m$ in the formula by $X$).
Note if $m=0$, all terms on the left of \eqref{combinatorics}
are zero but the first term $1$. Similarly, \eqref{formula1} is valid for all $m$, and
for negative $m$, $(1+x+x^2+x^3+\cdots )^{m} = (1-x)^{-m}$. Note that for $m=-1$,
all binomial coefficents $k\choose n$ appearing in the formula with $0\leq k<n$ are zero,
and we are left with $\displaystyle {-1\choose 1}  = -1$. This gives that
$(1+x+x^2+x^3+\cdots )^{-1}=1-x$, which is precisely what one expects.
Similarly the formula we use throughout
$\displaystyle {m\choose n-1} + {m\choose n} = {m+1\choose n}$ is valid for all $m$, and $n> 0$.
\ere

%\bre It might happen that $\brho$ is compact even if one of the weights $w_i>1$. For example if $w_1,\ldots, w_{r-1}\leq 1$, $[\rho]-1+w_r \leq\rho$, then any sequence $\sum t_ix_i$ with total weight $w$ converges to a sequence with total weight $w'$ that can be possibly be bigger than $w$ but never exceeds $\rho$. In this case as well, $\chi_c=\chi$. \ere

%%%%%%%%%%%%%%%%%%%%%%%%%%%%%%%%%%%%%%%%%%%%%%%%%%%%%%%%%%%%%%%%%%%%%%%%%%%%%%%%%%

\section{Barycenter Spaces of Disconnected Spaces}\label{disconnected}

We first observe that the Euler characteristic of ${\mathcal B}_n(X)$ only depends on $\chi = \chi (X)$, and not on the number of components.

\bpr \label{mcs} If $X$ has a finite number of components, then
$\displaystyle\chi {\mathcal B}_n(X) = 1- {n-\chi\choose n}$ for $n\geq 1$.
\epr

\begin{proof} For $n=1$ there is nothing to prove. Assume $n\geq 2$.
The proof proceeds by induction on the number of components. The formula is true for connected spaces \cite{kk}.
Write $X= Y\sqcup A$ where $A$ is a connected component. We can then
assume the theorem is true for $Y$. Pick $p_1\in Y$ and $p_2\in A$, and
consider the subspace
${\mathcal B}_{n}(X,\{p_1,p_2\})$ of all barycenter configurations containing either $p_1$ or $p_2$. This is then a union of two connected spaces
${\mathcal B}_{n}(X,\{p_1,p_2\}) = {\mathcal B}_{n}(X,p_1)\cup {\mathcal B}_{n}(X,p_2)$, and these spaces intersect
along the subspace
$${\mathcal B}_{n}(X)\cup {\mathcal B}_{n-1}(X,p_1,p_2)$$
(recall that ${\mathcal B}_0(X,p_1,p_2)={\mathcal B}_2(\{p_1,p_2\})$ which is an interval).
Since both spaces ${\mathcal B}_{n}(X,p_i)$ are contractible, their union has the homotopy type
of the suspension of their intersection and we have
\begin{eqnarray*}
{\mathcal B}_{n}(X,\{p_1,p_2\})&\simeq&\Sigma ({\mathcal B}_{n}(X)\cup {\mathcal B}_{n-1}(X,p_1,p_2)) \\
&\simeq&\Sigma \left({{\mathcal B}_{n}(X)\over {\mathcal B}_{n}(X)\cap {\mathcal B}_{n-1}(X,p_1,p_2)}\right)\ \ \ \ \
\ \hbox{since ${\mathcal B}_{n-1}(X,p_1,p_2)$ is contractible}\\
&\simeq&\Sigma \left({{\mathcal B}_{n}(X)\over {\mathcal B}_{n-1}(X,\{p_1,p_2\}}\right)\\
&\simeq&\Sigma {\mathcal B}_{n}(X/p_1\sim p_2)\ \ \ \ \ \hbox{by Lemma \ref{properties}, (ii)}\\
&\simeq&\Sigma {\mathcal B}_{n}(Y\vee A)
\end{eqnarray*}
Note that $Y\vee A$ has one component less than that of $X$ and $\chi (Y\vee A) = \chi - 1$, where
$\chi = \chi (X)$. Taking Euler characteristics gives that
\begin{eqnarray*}
\chi {\mathcal B}_{n}(X,\{p_1,p_2\}) &=& 2- \chi {\mathcal B}_{n}(Y\vee A) \\
&=&1+ {n - (\chi Y+\chi A-1)\choose n}\ \ \ \ \ \hbox{by induction}\\
&=&1 + {n - \chi+1 \choose n}
\end{eqnarray*}
On the other hand,
$\displaystyle
{{\mathcal B}_{n}(X)\over {\mathcal B}_{n-1}(X,\{p_1,p_2\})}\simeq {\mathcal B}_{n}(Y\vee A)$, and we can use
the formula for the quotient $\chi ({A\over B})=\chi (A) -\chi (B)+1$ to write
$$\chi {\mathcal B}_n(Y\vee A)=\chi \left({{\mathcal B}_{n}(X)\over {\mathcal B}_{n-1}(X,\{p_1,p_2\})}\right) =
\chi {\mathcal B}_n(X) +1 -\chi {\mathcal B}_{n-1}(X,\{p_1,p_2\})$$
which recombines into
\begin{eqnarray*}
\chi {\mathcal B}_n(X)&=& -1 +\chi {\mathcal B}_{n-1}(X,\{p_1,p_2\})
+\chi {\mathcal B}_n(Y\vee A)\\
&=& 1+ {n - \chi\choose n-1} - {n - \chi+1\choose n}=
1-{n-\chi\choose n}
\end{eqnarray*}
The proof is complete.
\end{proof}

\subsection*{A combinatorial proof of Proposition \ref{mcs}}
We know that the homology of ${\mathcal B}_k(X)$ only depends on the homology of $X$ (a fact more general than Euler characteristics) \cite{kk}. The following shows that this is true for disconnected spaces as well.

\begin{theorem}\label{cc}\cite{ahmedou}
Suppose $X=A\sqcup B$ is the disjoint union of spaces (not necessarily connected).
Then for $k\geq 2$, ${\mathcal B}_k(A\sqcup B)$ has the same homology as
\begin{eqnarray*}
&&{\mathcal B}_k(A)\vee \Sigma {\mathcal B}_{k-1}(A)\vee {\mathcal B}_k(B)\vee \Sigma {\mathcal B}_{k-1}B\\
&&\vee\bigvee_{\ell =1}^{k-1} {\mathcal B}_{k-\ell}(A)*{\mathcal B}_\ell (B)\ \vee\ \bigvee_{\ell =2}^{k-1}
\Sigma {\mathcal B}_{k-\ell}(A)*{\mathcal B}_{\ell -1}(B)
\end{eqnarray*}
\end{theorem}

\bex We can describe some homotopy types of some barycenter spaces of disjoint unions:
\begin{enumerate}[(i)]
\item When $k=2$, ${\mathcal B}_2(A\sqcup B)$ has actually the homotopy type of
${\mathcal B}_2(A)\vee \Sigma (A\times B)\vee {\mathcal B}_2(B)$.
\item When one of the components is contractible, say $B\simeq p$, then
$${\mathcal B}_n(A\sqcup B)\simeq {\mathcal B}_n(A\sqcup p)={\mathcal B}_n(A)\cup {\mathcal B}_{n-1}(A)*p
\simeq {\mathcal B}_n(A)\vee\Sigma {\mathcal B}_{n-1}A$$
This last equivalence follows from the fact that we are attaching a cone on ${\mathcal B}_{n-1}(A)$ which is itself contractible in ${\mathcal B}_n(A)$.
\item For $k\geq 2$, it is not hard to check that there is a homotopy equivalence
$${\mathcal B}_2(A\sqcup \{y_1,\ldots, y_k\})\simeq {\mathcal B}_2(A)\vee \bigvee^k\Sigma A\vee\bigvee^{k\choose 2}S^1
$$
which is in fact the decomposition in Theorem \ref{cc} obtained at the level of spaces
\item It is not always true that when the components are contractible, ${\mathcal B}_k(X)$ is also contractible for $k\geq 2$. This only happens when $k$ is larger (or equal) to the number of components. In fact, if $X=[n+1]$ is a set consisting of $n+1$-vertices, then ${\mathcal B}_{k+1}([n+1])$ is the $k$-th skeleton of $n$-dimensional simplex $\Delta_n$, and thus is a bouquet of spheres
    $${\mathcal B}_k([n+1])\simeq {{n\choose k}}S^{k-1}\ ,\ k\leq n$$
    which is the notation for a bouquet of that many spheres.
    This we can recover iteratively as follows (1)
    \begin{eqnarray*}
    {\mathcal B}_k([n+1])&=&{\mathcal B}_k([n]\sqcup [1]) \simeq {\mathcal B}_k([n])\vee \Sigma {\mathcal B}_{k-1}([n])\\
    &\simeq&\left({{n-1\choose k}}+ {{n-1\choose k-1}}\right) S^{k-1} ={{n\choose k}}S^{k-1}
    \end{eqnarray*}
\end{enumerate}
\eex

The combinatorial proof of Proposition \ref{mcs} now proceeds as follows.  By taking $\chi$ of the wedge decomposition in Theorem \ref{cc}, we find that
\begin{eqnarray*}
\chi {\mathcal B}_k(A\sqcup B) &=&\chi {\mathcal B}_kA + (2-\chi {\mathcal B}_{k-1}A) + \chi {\mathcal B}_kB + (2-\chi {\mathcal B}_{k-1}B)\\
&&+\sum_{\ell=1}^{k-1}(\chi {\mathcal B}_{k-\ell}A + \chi {\mathcal B}_\ell B - \chi {\mathcal B}_{k-\ell}A\chi {\mathcal B}_\ell B)\\
&&+\sum_{\ell=2}^{k-1}(2-\chi {\mathcal B}_{k-\ell}A - \chi {\mathcal B}_{\ell-1} B + \chi {\mathcal B}_{k-\ell}A\chi {\mathcal B}_{\ell-1} B)
-2k
\end{eqnarray*}
The term ``$-2k$" accounts for the wedge points being counted multiple times. We can set $\chi_1=\chi (A)$
and $\chi_2=\chi (B)$ and proceed by induction. Replacing binary coefficients in the above expression we get
the identity
\begin{eqnarray*}
&&\chi {\mathcal B}_k(A\sqcup B)\\
 &=& 4+ {k-1-\chi_1\choose k-1} - {k-\chi_1\choose k} + {k-1-\chi_2\choose k-1} - {k-\chi_2\choose k} \\
 &&+\sum_{\ell=1}^{k-1}\left(1-{k-\ell - \chi_1\choose k-\ell}{\ell-\chi_2\choose \ell}\right) + \sum_{\ell=2}^{k-1}\left(1 + {k-\ell-\chi_1\choose k-\ell}{\ell-1-\chi_2\choose\ell -1}\right) - 2k \\
&=&1 - {k-1-\chi_2\choose k} - {k-1-\chi_1\choose k}  - {k-1-\chi_1\choose k-1}{1-\chi_2\choose 1}\\
&&-\sum_{\ell=2}^{k-1}{k-\ell-\chi_1\choose k-\ell}\left[{\ell -\chi_2\choose\ell} -
{\ell -1 -\chi_2\choose \ell -1}\right]\\
&=&1 - {k-1-\chi_2\choose k} - {k-1-\chi_1\choose k}  - {k-1-\chi_1\choose k-1}{1-\chi_2\choose 1}
-\sum_{\ell=2}^{k-1}{k-\ell-\chi_1\choose k-\ell}{\ell -1 -\chi_2\choose\ell}\\
 &=&1-\sum_{\ell=0}^{k}{k-\ell - \chi_1\choose k-\ell}{\ell -1-\chi_2\choose \ell}\\
 &=&1-{k-\chi_1-\chi_2\choose k}
\end{eqnarray*}
The last identity can be verified directly, or found in \cite{gould}, formula (1.78). This completes this interesting computation and the proof of Proposition \ref{mcs}.

\subsection{Euler characteristic of $\brho$ for $X$ disconnected} We extend the previous computation to the
barycenter space with singular weights and to the Euler characteristic with compact supports. Here too we show that the final answer doesn't differ from the connected case. Starting with the stratification of $\brho$ in \eqref{disjointunion}, valid for all  $X$, we follow the steps in the proof of Proposition \ref{secondchic}. Only one step needs to be modified which is
the identification in Lemma \ref{identification} which is no longer true if $X$ is disconnected. The correct identification in that case is given by the following lemma.

\ble\label{quotient} Let $X$ be a disjoint union of components $A_i$, which are locally closed in $\overline{A}_i$, and of compact components ${\mathcal B}_j$.
We write $X=\bigsqcup_{i=1}^q A_i\sqcup\bigsqcup_{j=1}^t {B}_j$. Assume wlog that
$A_1,\ldots, A_s$ have singular points each of respective cardinality $a_1,\ldots, a_s\neq 0$. Assume as well that ${B}_1,\ldots, {B}_\ell$ have each singular points of respective cardinality $b_1,\ldots, b_\ell\neq 0$. Obviously $a_1 +\cdots +a_s+b_1 +\cdots +b_\ell=r$.
The other components $A_{s+1},\ldots, A_q, {B}_{\ell+1},\cdots , {B}_t$ have no singular points.
Then
$$\displaystyle {\overline{X}\over\partial\overline{X}\cup Q_r}\simeq
\displaystyle \left({\overline{A}_{1}\over\partial\overline{A}_1}\vee\cdots\vee {\overline{A}_{q}\over\partial\overline{A}_q}\vee {B}_1\vee\cdots {B}_\ell\vee \bigvee^{r-\ell} S^1\right)\bigsqcup
{B}_{\ell+1}\sqcup\cdots\sqcup {B}_{t}$$
\ele

\begin{proof} Each component ${B}_i$ with $b_i$ singular points contributes a bouquet $\bigvee^{b_i-1}S^1$ in the quotient.
Since $\partial\overline{X} = \bigsqcup\partial\overline{A_i}$,
each component $A_j$ with $a_j$ singular points contributes a bouquet of
$\bigvee^{a_i}S^1$ in the quotient (the extra leaf in the bouquet comes from the fact that there is non-trivial boundary, see Lemma \ref{identification}). Since $a_1+\cdots +a_s+b_1-1+\cdots b_\ell-1=r-\ell$, this accounts for the bouquet of circles and the rest is immediate.
\end{proof}

\bre\label{doublecheck} We can double-check this decomposition against the computation
$\chi_c(X-Q_r)=\chi_c(X)-r$. Indeed, let's recall that
$\chi (Z_1\vee\cdots\vee Z_n\vee\bigvee^m S^1)=\sum\chi (Z_i)-(n+m-1)$. We have
\begin{eqnarray*}
\chi \left({\overline{X}\over\partial\overline{X}\cup Q_r}\right)
&=&\sum_{i=1}^q\chi \left({\overline{A_i}\over \partial\overline{A_i}}\right) +
\sum_{j=1}^\ell\chi ({B}_i) - (q+\ell +r-\ell -1) + \sum_{j=\ell+1}^t\chi {B}_j \\
&=&\sum_{i=1}^q (\chi_c(A_i) + 1) + \sum_{j=1}^t\chi ({B}_j) -q -r  = \chi_c(X)-r+1
\end{eqnarray*}
where $\chi_c(X)=\sum \chi_c(A_i) + \sum\chi ({B})_j$.
This then gives that
$\chi_c(X-Q_r)=\chi \left({\overline{X}\over\partial\overline{X}\cup Q_r}\right)-1 = \chi_c(X)-r$
as expected.
\ere

\bth Theorem \ref{main} is true if $X$ is disconnected.
\end{theorem}

\begin{proof}
As in the proof of Proposition \ref{secondchic}, we use the same stratification
of $\brho$ with generic stratum
$X_n = {\mathcal B}_{n}(X-Q_r)\overset{\circ}{*}{\mathcal B}_k\{p_{i_1},\ldots ,p_{i_k}\}$.
By taking the appropriate compactification and applying both Lemmas \ref{mcs}, \ref{quotient}
and Remark \ref{doublecheck}, we obtain that
$$\chi_c(X_n) = (-1)^k\left( \chi {\mathcal B}_n\left( {\overline{X}\over \partial\overline{X}\cup Q_r}\right) -1\right)
= -(-1)^k{n-\chi_c + r-1\choose n}$$
which is the same as in the connected case. The rest of the argument runs
as in the proof in \S\ref{mainkey}.
\end{proof}

%%%%%%%%%%%%%%%%%%%%%%%%%%%%%%%%%%%%%%%%%%%%%%%%%%%%%%%%%%%%%%%%%%%%%%%%%%%%%%%%%%%%%%%%%%%%%%

\section{The Case of Two Singular Points} \label{degree2}

Interestingly, eventhough the distribution of the singular points among the components doesn't affect the Euler characteristic,  it does greatly affect the homology.
We analyze completely the various homotopy types (depending on weights) of the space ${\mathcal B}^{Q_2}_\rho (X)$ with $X$ having at most two connected components (this can easily be extended to any number of components).

\bpr\label{connectedr=2} Suppose $X$ is connected, $0<w_1\leq w_2\leq 1$.
Write $\rho = \lfloor\rho\rfloor +\epsilon$, $0\leq\epsilon<1$.
Then the possible homotopy types of ${\mathcal B}_\rho^{Q_2}(X)$ are :
\begin{enumerate}
\item contractible,  if $0 < w_1 + w_2\leq \epsilon$.
\item $\Sigma {\mathcal B}_{n}(X\vee S^1)$, if $w_1\leq\epsilon, w_2\leq \epsilon,
w_1 + w_2> \epsilon$.
\item contractible, if $w_1\leq \epsilon, w_2> \epsilon$.
\item ${\mathcal B}_{n}(X\vee S^1)$, if $w_1>\epsilon, w_2>\epsilon,
w_1 + w_2 \leq 1 + \epsilon$.
\item ${\mathcal B}_{n}(X)$, if $w_1+w_2 > 1+\epsilon$.
\end{enumerate}
\epr

\begin{proof}
Case (1): This is the case when
${\mathcal B}_\rho^{Q_2}(X)\simeq {\mathcal B}_{n}(X,p_1,p_2)$ which is contractible (Lemma \ref{contractible}).
Case (2): This is the case when
${\mathcal B}_\rho^{Q_2}(X) = {\mathcal B}_{n}(X,p_1)\cup {\mathcal B}_{n}(X,p_2)$.
This union was worked out in the proof of Proposition \ref{mcs} and is of the homotopy type of
$\Sigma {\mathcal B}_n(X/p_1\sim p_2)\simeq \Sigma {\mathcal B}_n(X\vee S^1)$.
Case (3): This is the case ${\mathcal B}_\rho^{Q_2}(X) = {\mathcal B}_{n}(X,p_1)$, and is contractible.
Case (4): This is the case ${\mathcal B}_\rho^{Q_2}(X) = {\mathcal B}_n(X)\cup {\mathcal B}_{n-1}(X,p_1,p_2)$
which is up to homotopy ${\mathcal B}_n(X\vee S^1)$.
Case (5): This is immediate.
\end{proof}

Turning to the disconnected case, we can write
$X = A_1\sqcup\cdots\sqcup A_q$ as a disjoint union of non-empty connected components. As before $Q_r = \{p_1,\ldots, p_r\}$ are the singular points with weights $w_1,\ldots, w_r$. Let's write
\begin{equation}\label{newnot}
{\mathcal B}_\rho^{Q_{r_1,\ldots ,r_q}}(X)
\end{equation}
for the subspace of $\brho$ consisting of configurations having $r_i$ of the singular
points in $A_i$, $\sum r_i=r$.

\bpr\label{disconnectedr=2} Let $X = A_1\sqcup A_2$, $r=2$. Then
\begin{enumerate}
\item Suppose $0 < \lfloor \rho\rfloor + w_1 + w_2\leq \rho$. Then
$$
{\mathcal B}_\rho^{Q_{1,1}}(X)\simeq {\mathcal B}_\rho^{Q_{2,0}}(X)\simeq {\mathcal B}_{\lfloor \rho\rfloor}(X,p_1,p_2)\simeq \ast
$$
\item $0 < \lfloor \rho\rfloor+w_1, \lfloor \rho\rfloor+w_2\leq \rho$,
$\lfloor \rho\rfloor+ w_1 + w_2> \rho$. Then
$$
{\mathcal B}_\rho^{Q_{1,1}}(X)\simeq \Sigma {\mathcal B}_{\lfloor \rho\rfloor}(A_1\vee A_2)$$
$${\mathcal B}_\rho^{Q_{2,0}}(X)\simeq \Sigma {\mathcal B}_{\lfloor \rho\rfloor}(A_1\vee S^1\sqcup A_2)$$
\item Suppose $0 < \lfloor \rho\rfloor+w_1\leq \rho$,
$\lfloor \rho\rfloor+w_2> \rho$. Then
$$
{\mathcal B}_\rho^{Q_{1,1}}(X)\simeq \ast\simeq {_1}B^{Q_{2,0}}(X)
$$
\item Suppose $\rho < \lfloor \rho\rfloor+w_1, \lfloor \rho\rfloor+w_2,
\lfloor \rho\rfloor+ w_1 + w_2 \leq 1 + \rho$. Then
$$
{\mathcal B}_\rho^{Q_{1,1}}(X)\simeq {\mathcal B}_{\lfloor \rho\rfloor}(A_1\vee A_2)$$
$${\mathcal B}_\rho^{Q_{2,0}}(X)\simeq {\mathcal B}_{\lfloor \rho\rfloor}(A_1\vee S^1\sqcup A_2)$$
\item Suppose $\rho < \lfloor \rho\rfloor+w_1, \lfloor \rho\rfloor+w_2\leq 1+\rho$, $\lfloor \rho\rfloor+w_1+w_2 > 1+\rho$. Then
$${\mathcal B}_\rho^{Q_{1,1}}(X)= {\mathcal B}_{\lfloor \rho\rfloor}(X)=
{\mathcal B}_\rho^{Q_{2,0}}(X)$$
\end{enumerate}
\epr

\begin{proof}
The proof runs exactly as in the previous proposition keeping track of
the identifications in ${\mathcal B}_n(X/p_1\sim p_2)$ as in Lemma \ref{quotient}.
\end{proof}

%%%%%%%%%%%%%%%%%%%%%%%%%%%%%%%%%%%%%%%%%%%%%%%%%%%%%%%
\vskip 20pt

%\addcontentsline{toc}{section}{Bibliography}
%\bibliography{biblio}

\begin{thebibliography}{doC}

\bibitem{ahmedou} M. Ahmedou, S. Kallel, C.B. Ndiaye, \textit{On the resonant boundary $Q$-curvature problem: the role of critical points at infinity and the boundary symmetric joins}, arXiv:1604.03745.
\bibitem{brown} R. Brown, \textit{Topology and Groupoids}, McGraw-Hill first version (1968). Online version www.groupoids.org.uk (2006).
\bibitem{carlotto} A. Carlotto, \textit{On the solvability of singular Liouville equations on compact surfaces of arbitrary genus}, Trans. Amer. Math. Soc. {\bf 366}, no.3, (2014), 1237--1256.
\bibitem{cm} A. Carlotto, A. Malchiodi, \textit{Weighted barycentric sets and singular Liouville equations on compact surfaces}, J. Funct. Anal. {\bf 262} (2012), no. 2, 409--450.
\bibitem{cl} C.C. Chen, C.S. Lin, \textit{Mean field equation of Liouville type with singular data: topological degree},  Comm. Pure Appl. Math. {\bf 68} (2015), no. 6, 887--947.
\bibitem{ginzburg} N. Chriss, V. Ginzburg, \textit{Representation theory and complex geometry}, Modern Birkh$\ddot{a}$user Classics. Birkh$\ddot{a}$user Boston, (2010).
\bibitem{dlr} F. De Marchis, R. L$\acute{o}$pez-Soriano, D. Ruiz, \textit{Compactness, existence and multiplicity for the singular mean field problem with sign-changing potentials}, J. Math. Pure Appl. {\bf 115} (2018), 237--267.
\bibitem{gr} M. Ganster, I.L Reilly, \textit{Locally closed sets and LC-continuous functions},
Internat. J. Math. and Math. Sci. {\bf 12}, 3 (1989), 417--424.
\bibitem{goresky} M.Goresky, \textit{Primer on sheaves}, lectures notes, IAS School of Mathematics.
\bibitem{gould} H.W. Gould, \textit{Combinatorial Identities: Table I: Intermediate
Techniques for Summing Finite Series}, on the internet (compiled by Jocelyn Quaintance).\
https://www.math.wvu.edu/$\sim$gould/Vol.4.PDF
\bibitem{kk} S. Kallel, R. Karoui, \textit{Symmetric Joins and Weighted Barycenters}, Adv. Nonlinear Stud. {\bf 11} (2011), 117--143.
\bibitem{kt} S. Kallel, W. Taamallah, \textit{On the Euler characteristic of stratified spaces}, work in progress.
\bibitem{kronenburg} M.J. Kronenburg, \textit{The binomial coefficient for negative arguments}, ArXiv 1105.3689.
\end{thebibliography}
\bibliographystyle{plain[8pt]}

\end{document}